\documentclass[11pt]{article}
\usepackage{cite}
\usepackage{mathrsfs}
\usepackage{amsfonts}
\usepackage{amsmath}
\usepackage{amsfonts,amssymb,color}
\usepackage{dsfont}
\usepackage{curves}
\usepackage{mathrsfs}
\usepackage{pifont}
\usepackage{amssymb}
\usepackage{latexsym,amsmath,amssymb,amsfonts,epsfig,graphicx,cite,psfrag}
\usepackage{eepic,color,colordvi,amscd}

\newtheorem{theorem}{Theorem}[section]

\newtheorem{lemma}{Lemma}[section]
\newtheorem{corollary}{Corollary}[section]

\newcommand{\qed}{\hfill\rule{0.5em}{0.809em}}

\def\emptyset{\mbox{{\rm \O}}}

\renewcommand{\baselinestretch}{1.2}
\def\qed{\hfill \rule{4pt}{7pt}}

\def\pf{\noindent {\it Proof. }}

\begin{document}
	
	\title{Structure and linear-Pollyanna for some square-free graphs}
	\author{Ran Chen\footnote{Email: 1918549795@qq.com},  \; Baogang  Xu\footnote{Email: baogxu@njnu.edu.cn. Supported by NSFC 11931006}\\\\
		\small Institute of Mathematics, School of Mathematical Sciences\\
		\small Nanjing Normal University, 1 Wenyuan Road,  Nanjing, 210023,  China}
	\date{}

\maketitle

\begin{abstract}
We use $P_t$ and $C_t$ to denote a path and a cycle on $t$ vertices, respectively. A {\em bull} is a graph consisting of a triangle with two disjoint pendant edges, a {\em hammer} is a graph obtained by identifying an endvertex of a $P_3$ with a vertex of a triangle. A class ${\cal F}$ is $\chi$-bounded if there is a function $f$ such that $\chi(G)\leq f(\omega(G))$ for all induced subgraphs $G$ of a graph in ${\cal F}$. A class ${\cal C}$ of graphs is  {\em Pollyanna} (resp. {\em linear-Pollyanna}) if ${\cal C}\cap {\cal F}$ is polynomially (resp. linearly) $\chi$-bounded for every $\chi$-bounded class ${\cal F}$ of graphs. Chudnovsky {\em et al} \cite{CCDO2023} showed that both the classes of bull-free graphs and hammer-free graphs are Pollyannas. Let $G$ be a connected graph with no clique cutsets and no universal vertices. In this paper, we show that $G$ is $(C_4$, hammer)-free if and only if it has girth at least 5, and $G$ is $(C_4$, bull)-free if and only if it is a clique blowup of some graph of girth at least 5. As a consequence, we  show that both the classes of $(C_4$, bull)-free graphs and $(C_4$, hammer)-free graphs are linear-Pollyannas. We also show that the class of (bull, diamond)-free graphs is linear-Pollyanna.

	\begin{flushleft}
	\end{flushleft}
	
\end{abstract}

\newpage

\section{Introduction}

All graphs considered in this paper are finite and simple. We follow \cite{BM08} for
undefined notations and terminologies. We use $N_G(v)$ to denote the set of vertices adjacent to $v$, $d_G(v)=|N_G(v)|$. Let $N_G(X)=\{u\in V(G)\setminus X\;|\; u$ has a
neighbor in $X\}$. If it does not cause any confusion, we usually omit the subscript $G$ and simply write $N(v)$, $d(v)$ and $N(X)$. For $X\subset V(G)$ and $v\in V(G)$, let $N_X(v)=N(v)\cap X$, and let $G[X]$ be the subgraph of $G$ induced by $X$.

We say that a graph $G$ contains a graph $H$ if $H$ is isomorphic to an induced subgraph of $G$. For a set ${\cal H}$ of graphs, we say that $G$ is ${\cal H}$-free if $G$ contains no $H\in {\cal H}$. If ${\cal H}=\{H_1,..., H_t\}$, we simply write $G$ is $(H_1,..., H_t)$-free instead.

For $u$, $v\in V(G)$, we write $u\sim v$ if $uv\in E(G)$, and write $u\not\sim v$ if $uv\not\in E(G)$. Let $v\in V(G)$, and let $X$ and $Y$ be two subsets of $V(G)$. We say that $v$ is {\em complete} to $X$ if $v$ is adjacent to all vertices of $X$, say that $v$ is {\em anticomplete} to $X$ if $v$ is adjacent to no vertex of $X$, and say that $X$ is complete (resp. anticomplete) to $Y$ if each vertex of $X$ is complete (resp. anticomplete) to $Y$.

A {\em clique} ({\em stable set}) of $G$ is a set of mutually adjacent (non-adjacent) vertices in $G$. The clique number of $G$, denoted by
$\omega(G)$, is the maximum size of a clique in $G$. A {\em clique cutset} of $G$ is a clique $K$ in $G$ such that
$G-K$ has more components than $G$, and a {\em cut-vertex} is a clique cutset of size 1.

{\em Substituting} a vertex $v$ of a graph $G$ by a graph $H$ is an operation which creates a new graph with vertex set $V(H)\cup V(G-v)$ and edge set $E(G-v)\cup \{xy~|~x\in V(H), y\in N_{G}(v)\}\cup E(H)$. When $H$ is a clique (not necessary nonempty), the substitution is said to be {\em blowing up} $v$ of $G$ into a clique. A graph $G$ obtained from a graph $H$ by blowing up all the vertices into cliques is said to be a {\em clique blowup } of $H$. Moreover, if each such a clique is nonempty, then it is said to be a {\em nonempty clique blowup} of $H$;
 {\em Under this literature, an induced subgraph of $H$ can be viewed as a clique blowup of $H$ with cliques of sizes 1 or 0}.

For a positive integer $k$, a $k$-{\em coloring } of $G$ is a function $\phi: V(G)\rightarrow \{1,\cdots,k\}$ such that $\phi(u)\ne \phi(v)$ if $u\sim v$. The {\em chromatic number} of $G$, denoted by $\chi(G)$, is the minimum number $k$ for which $G$ has a $k$-coloring.

A graph is {\em perfect} if all its induced subgraphs $H$ satisfy $\chi(H)=\omega(H)$. An induced cycle on at least 4 vertices is called  a {\em hole}, and its complement is called an {\em antihole}. A {\em $k$-hole} is a hole on $k$ vertices. A hole or antihole is {\em odd} or {\em even} if it has odd or even number of vertices. In 2006, Chudnovsky {\em et al} \cite{CRST2006} proved the Strong Perfect Graph Theorem, which confirms a long-standing conjecture of Berge \cite{CB1961}.
\begin{theorem}\label{perfect}{\em\cite{CRST2006}}
	A graph is perfect if and only if it is (odd hole, odd antihole)-free.
\end{theorem}

In 1975,  Gy\'{a}rf\'{a}s \cite{G75} proposed the following important concept. A class ${\cal F}$ is $\chi$-bounded if there is a function $f$ such that $\chi(G)\leq f(\omega(G))$ for all induced subgraphs $G$ of a graph in ${\cal F}$, then we call $f$ a {\em binding function} of ${\cal F}$.

A class of graphs is {\em polynomially} $\chi$-{\em bounded} if it has a polynomial binding function. There are many classes of graphs which are polynomially $\chi$-$bounded$ such as perfect graphs \cite{CRST2006}, even-hole free graphs \cite{CS2023}, and $T$-free graphs for some forest $T$ \cite{SR2019}. It is known that $P_3$-free graphs are union of complete graphs and so they are perfect, and $P_4$-free graphs are perfect \cite{S74}. But up to now, no polynomial binding function for $P_t$-free graphs has been found, when $t\geq 5$. The best known binding function for $P_5$-free graphs $G$ is $\omega(G)^{\log_2(\omega(G))}$ for $\omega(G)\ge3$ \cite{SSS2023}.

Esperet \cite{E2017} conjectured that every $\chi$-bounded class of graphs is
polynomially $\chi$-bounded. But this conjecture was recently disproved by
Bria\'{n}ski {\em et al} \cite{BDW2023}. Inspired by Esperet's conjecture,
Chudnovsky {\em et al} \cite{CCDO2023} considered its analog for proper
classes of graphs, and introduced the concepts of {\em Pollyanna} and
{\em strongly Pollyanna} as follows. A class ${\cal C}$ of graphs is {\em Pollyanna} if
${\cal C}\cap {\cal F}$ is polynomially $\chi$-bounded for every $\chi$-bounded class
${\cal F}$ of graphs. For a positive integer $n$, a class ${\cal F}$ of graphs is $n$-{\em good} if it is hereditary and there is some constant $m$ such that $\chi(G)\leq m$ for every graph $G\in{\cal F}$ with $\omega(G)\leq n$. A class ${\cal C}$ of graphs is {\em $n$-strongly Pollyanna} if ${\cal C}\cap {\cal F}$ is polynomially $\chi$-bounded for every $n$-good class ${\cal F}$ of graphs. And, ${\cal C}$ is {\em strongly Pollyanna} if it is $n$-strongly Pollyanna for some integer $n$. It is clear that strongly Pollyannas are all Pollyannas.

\begin{figure}[htbp]
	\begin{center}
		\includegraphics[width=15cm]{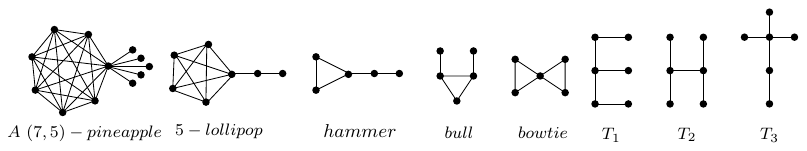}
	\end{center}
	\vskip -25pt
	\caption{Illustration of some special forbidden graphs.}
\label{fig-1}
\end{figure}

For positive integers $t$ and $k$, a {\em $(t,k)$-pineapple} is a graph obtained by
attaching $k$ pendant edges to a vertex of a $K_t$, a {\em $t$-lollipop} is a graph obtained from the disjoint union of a $K_t$ and a $P_2$ by adding an edge with one end in $K_t$ and the other end in $P_2$, a {\em hammer} is a 3-lollipop, a {\em bowtie} is a graph obtained from two copies of $K_3$ by sharing exactly one vertex, a {\em bull} is the graph consisting of a $K_3$ with two disjoint pendant edges (see Figure~\ref{fig-1}).

\renewcommand{\baselinestretch}{1}
\begin{theorem}\label{Chud-pollyanna}{\em\cite{CCDO2023}}
Let $m,k,t$ be positive integers.  Then,
\begin{itemize}
\item the class of $mK_t$-free graphs is $(t-1)$-strongly Pollyanna,

\item the class of $(t,k)$-pineapple-free graphs is $(2t-4)$-strongly Pollyanna,

\item the class of $t$-lollipop-free graphs is $(3t-6)$-strongly Pollyanna,

\item the class of bowtie-free graphs is $3$-strongly Pollyanna, and

\item the class of bull-free graphs is $2$-strongly Pollyanna.
\end{itemize}
\end{theorem}\renewcommand{\baselinestretch}{1.2}

A graph $H$ is {\em Pollyanna}-{\em binding} if the class of $H$-free graphs is Pollyanna. To investigate the structure of pollyanna-binding graphs, Chudnovsky {\em et al} \cite{CCDO2023} introduced the concept of {\em willow} (see \cite{CCDO2023} for its definition). Among other open questions on Pollyannas, they posed a very beautiful conjecture: {\em a graph is pollyanna-binding if and only if it is a willow}.

A class of graphs is {\em linearly} $\chi$-{\em bounded} if it has a linearly polynomial binding function. Both perfect graphs \cite{CRST2006} and even-hole free graphs \cite{CS2023} are linearly $\chi$-bounded.

There are so many classes of graphs which are polynomially $\chi$-bounded, such as $(P_6,C_4)$-free graphs \cite{GH2019}, $(P_6$, diamond)-free graphs \cite{CHM2021} (diamond is a graph obtained from a $P_3$ by adding a new vertex adjacent to all the vertices of $P_3$) and (fork, $C_4$)-free graphs \cite{CHKK}. What is the optimal binding function of a polynomially $\chi$-bounded class of graphs$?$ It is interesting to know whether a polynomially $\chi$-bounded
class has a linearly binding function or not.

We say that a class ${\cal C}$ of graphs is {\em linear-Pollyanna} if ${\cal C}\cap {\cal F}$ is linearly $\chi$-bounded for every $\chi$-bounded class ${\cal F}$ of graphs. Let $n$ be a positive integer. We say that a class ${\cal C}$ of graphs is {\em $n$-strongly linear-Pollyanna} if ${\cal C}\cap {\cal F}$ is linearly $\chi$-bounded for every $n$-good class ${\cal F}$ of graphs. And, we say that ${\cal C}$ is {\em strongly linear-Pollyanna} if it is $n$-strongly linear-Pollyanna for some integer $n$.

Let $P_k$ and $C_k$ denote a path and a cycle on $k$ vertices, respectively. In this paper, we prove that, for the graph $H$ which is a bull or a hammer, the class of $(C_4, H)$-free graphs is  2-strongly linear-Pollyanna. So, both classes are linear-Pollyanna.

\begin{theorem}\label{pollyanna}
Let $H$ be a bull or a hammer, and let ${\cal C}$ be the class of $(C_4, H)$-free graphs. Then ${\cal C}$ is $2$-strongly linear-Pollyanna.
\end{theorem}

Theorem~\ref{pollyanna} is a direct consequence of the following stronger conclusion.

\begin{theorem}\label{bull,hammer-pollyanna}
Let $H$ be a bull or a hammer, and ${\cal C}$ be the class of $(C_4, H)$-free graphs. Let $k\geq 2$ be an
integer, and ${\cal F}$ be a hereditary class of graphs satisfying $\chi(F)\leq k$ for every $F\in {\cal F}$ with $\omega(F)\leq 2$. Then, for every $G\in {\cal C}\cap {\cal F}$, $\chi(G)\leq \frac{k}{2}\omega(G)$ if $H$ is a bull, and $\chi(G)\leq\omega(G)+k-2$ if $H$ is a hammer.
\end{theorem}

To prove Theorem~\ref{bull,hammer-pollyanna}, we will first study the structure of $(C_4, H)$-free graphs for $H$ being a bull or a hammer, and prove the following theorem. A vertex $v$ of a graph $G$ is called a {\em universal vertex} if $v$ is complete to $V(G)\setminus \{v\}$. The {\em girth} of $G$ is the length of a shortest cycle in $G$.

\begin{theorem}\label{bull,hammer-1}
Let $G$ be a connected graph which has no clique cutsets or universal vertices. Then, $G$ is $(C_4,$ hammer)-free if and only if it has girth at least 5, and $G$ is $(C_4,$ bull)-free if and only if it is a clique blowup of some graph of girth at least $5$.
\end{theorem}

Let $F$ be a graph, and let ${\cal F}$ be the class of $(F, C_4)$-free graphs. By taking ${\cal C}$ to be the $(C_4, bull)$-free graphs and $(C_4, hammer)$-free graphs respectively, and applying Theorem~\ref{bull,hammer-pollyanna} to ${\cal C}$ and ${\cal F}$, we can immediately obtain the following theorem.

\begin{theorem}\label{bull,hammer-2}
Let $F$ be a graph and $k\geq 2$ be an
integer. If $(F, C_4, C_3)$-free graph is $k$-colorable, then $\chi(G)\leq \frac{k}{2}\omega(G)$ whenever $G$ is $(F, C_4, bull)$-free, and $\chi(G)\leq\omega(G)+k-2$ whenever  $G$ is $(F, C_4, hammer)$-free.
\end{theorem}

Let $T_1$ be the graph obtained from $K_{1,3}$ by subdividing two edges, $T_2$ be
the connected graph with two vertices of degree 3 and four vertices of degree of 1, and
$T_3$ be the graph obtained from $K_{1,4}$ by subdividing one edge (see Figure 1).
Randerath \cite{R1998} (see \cite{RS2004} Theorem~43) proved that for each $i\in\{1, 2, 3\}$,  $(T_i, C_3)$-free graphs
are 3-colorable, and Chudnovsky {\em et al} \cite{MCJS2018} proved that $(P_8, C_4, C_3)$-free graphs are 3-colorable. Apply Theorem~\ref{bull,hammer-2} by taking $F\in\{P_8, T_1, T_2, T_3\}$, we can conclude the following corollary.

\begin{corollary}\label{tree}
Let $G$ be an $(F, C_4, H)$-free graph for some $F\in\{P_8, T_1, T_2, T_3\}$ and
$H\in\{bull,hammer\}$. Then, $\chi(G)\leq\frac{3}{2}\omega(G)$ if $H$ is a bull, and $\chi(G)\leq\omega(G)+1$ if $H$ is a hammer.
\end{corollary}

In \cite{GH2019}, Gaspers and Huang proved that every $(P_6, C_4)$-free graph $G$ satisfies $\chi(G)\le {3\omega(G)\over 2}$ (Karthick and Maffray improved this bound to $\lceil\frac{5}{4}\omega(G)\rceil$ \cite{KM2019}), which indicates that $(P_6, C_4, C_3)$-free graphs are 3-colorable. Apply Theorem~\ref{bull,hammer-2} by taking $F=P_6$, we have

\begin{corollary}
Every $(P_6, C_4$, hammer)-free graphs is $(\omega(G)+1)$-colorable.
\end{corollary}

Let $T$ be a tree on $k$ vertices. In \cite{gyarfas2}, Gy\'{a}rf\'{a}s {\em et al} proved that every $(T, C_4, C_3)$-free graph has chromatic number at most $k-1$. As another direct consequence of Theorem~\ref{bull,hammer-2}, we have

\begin{corollary}\label{coro-tree-1}
Let $k\ge 3$, $T$ be a tree on $k$ vertices, and let $G$ be a $(T, C_4, H)$-free graph with $H$ being a bull or a hammer. Then, $\chi(G)\leq\frac{k-1}{2}\omega(G)$ if $H$ is a bull, and $\chi(G)\leq\omega(G)+k-3$ if $H$ is a hammer.
\end{corollary}

We will prove Theorems~\ref{bull,hammer-pollyanna} and \ref{bull,hammer-1} in the following Section 2. Section 3 is devoted to (bull, diamond)-free graphs, in which we show that the class of (bull, diamond)-free graphs is also linear-Pollyanna.

\section{Proofs of Theorems~\ref{bull,hammer-pollyanna} and \ref{bull,hammer-1}}

In this section, we prove Theorems~\ref{bull,hammer-pollyanna} and \ref{bull,hammer-1}.

Let $G$ be a graph and $F$ be an induced subgraph of $G$. We call $v$ a {\em center} for $F$ if $v\notin V(F)$ and $v$ is complete to $V(F)$, and call $v$ a {\em simplicial vertex} of $G$ if $N_G(v)$ is a clique. For two vertices $u$ and $v$, let $d(u, v)$ be the distance of $u$ and $v$, which is the length of a shortest path between $u$ and $v$. Let $S_1,S_2\subseteq V(G)$. We define $d(S_1,S_2)=\min\{d(s_1,s_2)~|~s_1\in S_1,s_2\in S_2\}$.

We call a subset $S$ of vertices a {\em homogeneous set} of $G$ if $2\leq|S|\leq |V(G)|-1$ and every vertex of $V(G)\setminus S$ is either complete  or anticomplete to $S$. A {\em homogeneous clique} is a homogeneous set which is a clique. A {\em chordal graph} is one without holes.

Let $H$ be a bull or a hammer, and let ${\cal C}$ be the class of $(C_4, H)$-free graphs. We begin from graphs of ${\cal C}$ with some further restrictions.

\begin{lemma}\label{C_4,bull,hammer}
Let $H$ be a bull or a hammer, and let $G$ be a connected $(C_4,H)$-free graph. If $G$ contains a hole with a center, then  $G$ has either a clique cutset or a universal vertex.
\end{lemma}
\pf  Suppose that $C=v_1v_2\cdots v_q v_1$ is a hole of $G$ and $v$ is a center for $C$, where $q\ge 5$. Let $I=\{x\in V(G)~|~x~\mbox{is a center for}~C\}$. Clearly, $I\ne\emptyset$ because $v\in I$, and $I$ must be a clique since any pair $x$ and $x'$ of nonadjacent vertices in $I$ would produce a 4-hole $xv_1x'v_3x$.

Let $G'=G-I$. For each nonnegative integer $i$, let $N_i=\{x\in V(G')~|~d(x,V(C))=i\}$.
Let $X=\bigcup_{i\geq0}N_i$, and $Y=V(G')\setminus X$. It is certain that $I$ is a clique
cutset of $G$ if $Y\ne\emptyset$ because $G$ is connected. Therefore, we may assume
that $Y=\emptyset$, and now $V(G')=X=\bigcup_{i\geq0}N_i$. Obviously, $N_0=V(C)$ is
complete to $I$. We will prove that every vertex in $I$ is a universal vertex of $G$.

First, we prove that
\begin{equation}\label{N_1}
	\mbox{$N_1$ is complete to $I$.}
\end{equation}

Let $x\in N_1$. Without loss of generality, we may assume that $x\sim v_1$. Suppose that $x$ has a non-neighbor, say $y$, in $I$. For each $i_{0}\in \{3, 4, \ldots, q-1\}$, to forbid a 4-hole $xv_{i_0}yv_1x$, we have that $x\not\sim v_{i_0}$. Therefore, $x$ is anticomplete to $V(C)\setminus\{v_1,v_2,v_q\}$. If $x\sim v_2$ and $x\sim v_{q}$, then a 4-hole $xv_qyv_2x$ appears. By symmetry, we may  assume that $x\not\sim v_2$. But then $G[\{v_1,v_2,v_4,x,y\}]$ is a bull, and $G[\{v_1,v_3,v_4,x,y\}]$ is a hammer, a contradiction. This proves (\ref{N_1}).

Next, we show that
\begin{equation}\label{N_i}
	\mbox{$N_i$ is complete to $I$ for $i\geq 1$.}
\end{equation}

Suppose to its contrary, we choose $i$ to be the minimum integer such that $N_i$ has a vertex $x$ that is not complete to $I$, and let $y\in I$ be a non-neighbor of $x$. By (\ref{N_1}), $i\geq2$. Let $P=x_ix_{i-1}x_{i-2}\cdots x_1x_0$ be a shortest path from $x$ to $V(C)$, where $x_0\in V(C)$ and $x_i=x$. We may assume that $x_0=v_1$. 

If $x_1\not\sim v_2$, then $x_1\not\sim v_3$ to avoid a 4-hole $x_1v_1v_2v_3x_1$, and so $G[\{x_i,x_{i-1},x_{i-2},v_3,y\}]$ is a bull and $G[\{x_i,x_{i-1},v_2,v_3,y\}]$ is a hammer, a contradiction. Therefore, $x_1\sim v_2$, and $x_1\sim v_{q}$ by symmetry.

Since $x_1\notin I$, it has a non-neighbor in $C$. Let $t\in \{3,4,\cdots,q-1\}$ be the smallest integer such that $x_1\not\sim v_t$. It is certain that $t\leq q-2$ to avoid a 4-hole $x_1v_{q-2}v_{q-1}v_{q}x_1$, and $x_1\not\sim v_{t+1}$ to avoid a 4-hole $x_1v_{t-1}v_{t}v_{t+1}x_1$. Then, $G[\{x_i, x_{i-1},x_{i-2},v_t,y\}]$ is a bull and $G[\{x_1,v_{t-2},v_{t-1},v_{t},v_{t+1}\}]$ is a hammer, a contradiction. This proves (\ref{N_i}).

It follows directly from (\ref{N_1}) and (\ref{N_i}) that $V(G')$ is complete to $I$, and thus every vertex in $I$ is a universal vertex of $G$. This completes the proof of Lemma~\ref{C_4,bull,hammer}. \qed

\begin{lemma}\label{C_4,bull}
Let $G$ be a connected $(C_4,bull)$-free graph without clique cutsets or universal vertices, $x\in V(G)$, and $C$ be a hole of $G-x$. If $x$ has two consecutive neighbors on $C$, then $G$ has a homogeneous clique.
\end{lemma}
\pf Suppose that $C=v_1v_2\cdots v_q v_1$ for some $q\ge 5$, and suppose that $x\sim v_1$ and $x\sim v_2$. Then, $x$ must have a neighbor in $\{v_3, v_q\}$ to avoid an induced bull on $\{x, v_1, v_2, v_3, v_q\}$. Suppose by symmetry that $x\sim v_3$.

By Lemma~\ref{C_4,bull,hammer}, $C$ has no center since $G$ has no clique cutsets and universal vertices. Thus, $x$ is not complete to $V(C)$. Let $j$ be the smallest integer in $\{4,5,\cdots,q\}$ such that $v_j\not\sim x$. If $j\geq5$, then $G[\{x,v_1,v_{j-2},v_{j-1},v_{j}\}]$ is a bull if $j\neq q$, or $xv_{j-1}v_jv_1x$ is a 4-hole if $j=q$.
Therefore, $j=4$, and $x\not\sim v_5$ to avoid a 4-hole $xv_3v_4v_5x$. If $x$ has a neighbor in $V(C)\setminus \{v_1,v_2,v_3\}$, then $q\geq 6$,
and suppose $x\sim v_k$ for some $k\in \{6, \ldots, q\}$; but then $G[\{x, v_2, v_3, v_4, v_k\}]$ is
a bull. Therefore, we have that
$$\mbox{$N_{V(C)}(x)=\{v_1, v_2, v_3\}$.}$$

Next, we show that $\{x,v_2\}$ is a homogeneous clique of $G$. Suppose to its contrary
that there is a vertex $y\in V(G)\setminus\{x, v_2\}$ such that $y$ has exactly one neighbor
in $\{x,v_2\}$. By symmetry, we may assume that $y\sim x$ and $y\not\sim v_2$.  If
$y\not\sim v_3$, then $G[\{x, y, v_2, v_3, v_4\}]$ is a bull whenever $y\not\sim v_4$,
and $xyv_4v_3x$ is a 4-hole whenever $y\sim v_4$. Therefore, $y\sim v_3$, and $y\sim v_1$ by
symmetry. But then, a 4-hole $v_1v_2v_3yv_1$ appears. This shows that every vertex of
$V(G)\setminus\{x, v_2\}$ is either complete or anticomplete to $\{x, v_2\}$, and thus
Lemma~\ref{C_4,bull} holds. \qed

\begin{lemma}\label{triangle}
Let $G$ be a connected graph which contains a triangle and a hole. If $G$ is 2-connected, then there exist a vertex $x\in V(G)$ and a hole $C$ such that $x$ has two consecutive neighbors on $C$.
\end{lemma}
\pf Suppose that $G$ is 2-connected. We choose a triangle $H_1=xyzx$ and a hole  $H_2=v_1v_2\cdots v_q v_1$  such that $d(V(H_1), V(H_2))$ is minimum among all triangles and holes of $G$.
Let $d(V(H_1), V(H_2))=t$. We first prove that
\begin{equation}\label{distance}
	t=0.
\end{equation}

Suppose to its contrary that $t\geq1$. Let $P=x_1x_2\cdots x_{t+1}$ be a shortest path from $H_1$ to $H_2$, and by symmetry suppose that $x_1=x$ and $x_{t+1}=v_1$, and this implies that $d(x,H_2)=d(H_1,H_2)\leq \min\{d(y,H_2),d(z,H_2)\}$. Hence, $\{z,y\}$ is anticomplete to $V(P)\setminus\{x_1,x_2\}$. If $\{z,y\}$ is not anticomplete to $x_2$, then by symmetry suppose $y\sim x_2$; but now $d(H_1,G[\{x,y,x_2\}])<d(H_1,H_2)$ which contradicts with the choice of $H_1$ and $H_2$. Therefore, $\{y, z\}$ must be anticomplete to $V(P)\setminus\{x_1\}$.

Since $G-x_1$ is connected as $G$ is 2-connected, we may choose a shortest path in $G-x_1$ between $\{y, z\}$ and $\{x_2\}$, say $P'=y_1y_2\cdots y_{t'}$. Suppose that $y_{t'}=x_2$, and  $y_1=y$ by symmetry. It is certain that $t'\ge 3$ as $\{y, z\}$ is anticomplete to $V(P)\setminus\{x_1\}$. If $x_1\sim y_{t'-1}$, then $x_1x_2y_{t'-1}x_1$ is a triangle that is closer to $H_2$ than $H_1$, contradicting the choice of $H_1$. Therefore, $x_1\not\sim y_{t'-1}$. Notice that $x_1=x$ and $y_1=y$. Let $i$ be the maximum integer in $\{1,2,\cdots,t'-2\}$ such that $x_1\sim y_{i}$. Then, $x_1y_iy_{i+1}\cdots y_{t'}x_1$ is a hole of $G$ that has a common vertex $x$ with $H_1$, contradicting $t\ge 1$. This shows that $t=0$, and so (\ref{distance}) holds.

We may by symmetry suppose that $x=v_1$. If $y,z\in V(H_2)$ or $\{y, z\}$ is not anticomplete to $\{v_2, v_q\}$, then we are done. So, we may suppose that $V(H_1)\cap V(H_2)=\{x\}$, and $\{y, z\}$ is anticomplete to $\{v_2, v_q\}$.

Since $G-x$ is connected, we may choose  $P''=z_1z_2\cdots z_{\ell}$ to be a shortest path
in $G-x$ from $\{y, z\}$ to $\{v_2\}$. It is certain that $\ell\ge 3$. Suppose that
$z_{\ell}=v_2$, and $z_1=z$ by symmetry. Without loss of generality, we may assume that
$x\not\sim z_{2}$, as otherwise we can replace $H_1$ by $xz_1z_2x$ and repeat the procedure.
Let $j$ be the smallest integer in $\{3,4,\cdots,{\ell}\}$ such that $x\sim z_{j}$. Then, $xz_1z_2\cdots z_{j}x$ is a hole such that $y\sim x$ and $y\sim z_1$. This proves Lemma~\ref{triangle}. \qed

\medskip

To prove Theorem~\ref{bull,hammer-1}, we also need the following lemma.

\begin{lemma}\label{simplicial}{\em\cite{D1961}}
Each chordal graph has a simplicial vertex.
\end{lemma}

\medskip

\noindent\textbf{{\em Proof of Theorem~\ref{bull,hammer-1}}:} The sufficiency is trivial. We need only to prove its necessity. Let $H$ be a bull or a hammer, and let $G$ be a connected $(C_4, H)$-free graph without clique cutsets or universal vertices. We will show that $G$ has girth at least 5 if $H$ is a hammer, and $G$ is a clique blowup of some graph of girth at least 5 if $H$ is a bull.

If $G$ does not contain a hole, then $G$ is a chordal graph, and thus has a simplicial vertex, say $v$, by Lemma~\ref{simplicial}. Since $N_G(v)$ is a clique, we have that $\{v\}$ is a universal vertex of $G$ if $V(G)=N_G(v)\cup\{v\}$, and $N_G(v)$ is a clique cutset of $G$ otherwise. Both are contradictions. Therefore, $G$ contains a hole. We now divide the proof into two cases depending on $H$ is a bull or a hammer.

\medskip

\noindent{\bf Case 1} First we suppose that $H$ is a hammer.

We need only to prove that $G$ is triangle-free. Suppose to its contrary that $G$ contains a triangle. By Lemma~\ref{triangle}, $G$ has a vertex $x$ and a hole $C=v_1v_2\cdots v_qv_1$ such that $x$ has two consecutive neighbors on $C$, where $q\geq 5$. Suppose by symmetry that $x\sim v_1$ and $x\sim v_2$.

We can deduce that $x\sim v_3$ as otherwise either a 4-hole $xv_2v_3v_4x$ (when $x\sim v_4$) or  a hammer on $\{x,v_1,v_2,v_3,v_4\}$ (when $x\not\sim v_4$) appears. Similarly, $x\sim v_q$. Since $G$ has no clique cutsets and universal vertices,  it follows from Lemma~\ref{C_4,bull,hammer} that $x$ cannot be complete to $V(C)$. Let $k$ be the smallest integer in $\{4,5,\cdots,q-1\}$ such that $x\not\sim v_k$. Then, $G[\{x,v_{k-2},v_{k-1},v_{k},v_{k+1}\}]$ is a hammer if $x\not\sim v_{k+1}$, and $xv_{k-1}v_kv_{k+1}x$ is a 4-hole if $x\sim v_{k+1}$. Both are contradictions. Therefore, $G$ has girth at least 5.

\medskip

\noindent{\bf Case 2} Now, we suppose that $H$ is a bull.

We prove this statement by induction on $|V(G)|$. If $|V(G)|$ is small enough, then $G$ has girth at least 5, and so the lemma holds.  Now, we proceed with the inductive proof.

If $G$ has no homogeneous cliques, then  $G$ is triangle-free by Lemmas~\ref{C_4,bull} and \ref{triangle}, and so $G$ has girth at least 5. Suppose that $G$ has a homogeneous clique, and let $S$ be a maximal homogeneous clique of $G$. Let $G'$ be the graph obtained from $G-S$ by adding a new vertex $v$ complete to $N_G(S)$ and anticomplete to $V(G)\setminus N_G(S)$. We have that $G$ can obtained from $G'$ by substituting $v$ of $G'$ by a clique $S$. Therefore, $G'$ has no clique cutsets and no universal vertices as $G$ has these properties. Also, $G'$ is a connected $(C_4$, bull)-free graph as $G'$ is isomorphic to some induced subgraph of $G$. Hence, by the inductive hypothesis, $G'$ is a clique blowup of some graph of girth at least 5. There exists a connected $(C_3,C_4)$-free graph $H$ such that $G'$ is a nonempty clique blowup of $H$. Since $H$ is isomorphic to some induced subgraph of $G'$, we may assume that $H$ is an induced subgraph of $G'$.

\begin{figure}[htbp]
	\begin{center}
		\includegraphics[width=8cm]{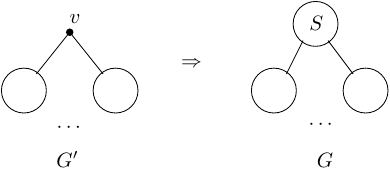}
	\end{center}
	\vskip -25pt
	\caption{Illustration of $G'$ and $G$.}
	\label{fig-2}
\end{figure}

By the maximality of $S$, we have that $v$ does not belong to any homogeneous clique of $G'$, and so $v\in V(H)$. Recall that $G$ can obtained from $G'$ by substituting $v$ of $G'$ by a clique $S$, and so $G$ is a clique blowup of $H$ (see Figure~\ref{fig-2}). This completes the proof of Theorem~\ref{bull,hammer-1}. \qed

\medskip

Before proving Theorem~\ref{bull,hammer-pollyanna}, we need the following lemma.

\begin{lemma}\label{clique blowup}
Let $k\ge 2$, and let $H$ be a graph with $\chi(H)\leq k$. If $G$ is a clique blowup of $H$,
then $\chi(G)\leq \frac{k}{2}\omega(G)$.
\end{lemma}
\pf We prove the conclusion by induction on $|V(G)|$. If $G$ is an induced subgraph of $H$ (this means
that $G$ is a clique blowup of $H$ with cliques of size 1 or 0),
then $\chi(G)\leq k$ whenever $\omega(G)\ge 2$, and $\chi(G)=1\le \frac{k}{2}$ whenever
$\omega(G)=1$. Now we proceed to inductive proof, and suppose, without loss of generality,
that $G$ is connected. Let $H'$ be an induced subgraph of $H$ such that $G$ is a nonempty
clique blowup of $H'$. Obviously, $H'$ is connected and $\chi(H')\leq k$.
If $\omega(H')\geq2$, then we have that $\omega(G-V(H'))\leq\omega(G)-2$, and thus by
inductive hypothesis, $\chi(G)\leq\chi(G-V(H'))+\chi(H')\leq\frac{k}{2}\omega(G-V(H'))+k\leq\frac{k}{2}(\omega(G)-2)+k=\frac{k}{2}\omega(G)$. If $\omega(H')=1$, then $G$ is a complete graph, and thus $\chi(G)=\omega(G)\leq \frac{k}{2}\omega(G)$. This proves Lemma~\ref{clique blowup}. \qed

Now, we are ready to prove Theorem~\ref{bull,hammer-pollyanna}.

\medskip

\noindent\textbf{{\em Proof of Theorem~\ref{bull,hammer-pollyanna}}: } Let $H$ be a bull or a hammer, and $k\ge 2$ be an integer. Let ${\cal C}$ be the class of $(C_4, H)$-free graphs, and let ${\cal F}$ be a hereditary class of graphs satisfying $\chi(F)\leq k$ whenever $\omega(F)\leq 2$ for every $F\in {\cal F}$.  Let $G\in {\cal C}\cap {\cal F}$ be a graph. We will prove the theorem by induction on $|V(G)|$. The conclusion is certainly true when $G$ has a small number of vertices.

If $G$ is disconnected, then there exists a component $G'$ of $G$ such that $\chi(G)=\chi(G')$. By inductive hypothesis, the Theorem holds.

If $G$ has a clique cutset $Q$, let $V(G)\setminus Q$ be partitioned into two nonempty subsets $A$ and $B$ such that $\chi(G[Q\cup A])\ge \chi(G[Q\cup B])$, then $\chi(G)=\chi(G[Q\cup A])\leq \frac{k}{2}\omega(G[Q\cup A])\leq\frac{k}{2}\omega(G)$ if $H$=bull, and $\chi(G)=\chi(G[Q\cup A])\leq\omega(G[Q\cup A])+k-2\leq\omega(G)+k-2$ if $H$=hammer.

If $G$ has a universal vertex $v$, then $\omega(G)=\omega(G-v)+1$. By inductive hypothesis, $\chi(G)=\chi(G-v)+1\leq \frac{k}{2}\omega(G-v)+1= \frac{k}{2}(\omega(G)-1)+1\leq\frac{k}{2}\omega(G)$ if $H$=bull, and $\chi(G)=\chi(G-v)+1\leq\omega(G-v)+k-2+1=\omega(G)+k-2$ if $H$=hammer.

So, we suppose that $G$ is a connected graph which has no clique cutsets and universal vertices.

If $H$ is a bull, then by Theorem~\ref{bull,hammer-1}, $G$ is a clique blowup of a graph $F$ of girth at least 5. Since $F$ is an induced subgraph of $G$, we have that $F\in {\cal F}$, and so $\chi(F)\leq k$. By Lemma~\ref{clique blowup}, $\chi(G)\leq\frac{k}{2}\omega(G)$.

If $H$ is a hammer, then by Theorem~\ref{bull,hammer-1}, $G$ has girth at least 5, and thus $\chi(G)\leq k$. It is obviously that $\chi(G)\leq k= \omega(G)+k-2$ if $\omega(G)=2$, and  $\chi(G)=1\leq\omega(G)+k-2$ if $\omega(G)=1$. This completes the proof of Theorem~\ref{bull,hammer-pollyanna}. \qed

\section{The class of (bull, diamond)-free graphs}

In this section, we will show that the class of (diamond, bull)-free graphs is 2-strongly linear-Pollyanna.

\begin{theorem}\label{diamond-1}
	Let ${\cal C}$ be the class of (bull, diamond)-free graphs. Let $k\geq 2$ be an
	integer, and ${\cal F}$ be a hereditary class of graphs satisfying $\chi(F)\leq k$ for every $F\in {\cal F}$ with $\omega(F)\leq 2$. Then, for every $G\in {\cal C}\cap {\cal F}$, $\chi(G)\leq \max\{k,\omega(G)\}$, and the class ${\cal C}$ is $2$-strongly linear-Pollyanna.
\end{theorem}

As usual, we use $\delta(G)$ to denote the minimum degree of $G$. The {\em Cartesian product} of any two graphs $G$ and $H$, denoted by $G\square H$, is the graph with vertex set $\{(a,u)~|~a\in V(G)~\mbox{and}~b\in V(H)\}$, where two vertices $(a,u)$ and $(b,v)$ are adjacent if either $a=b$ and $u\sim v$ in $H$, or $u=v$ and $a\sim b$ in $G$. To prove Theorem~\ref{diamond-1}, the following theorem is very useful.

\begin{theorem}\label{diamond-2}
	Let $G$ be a connected (bull, diamond)-free graph. If $\omega(G)\geq3$, then either $\delta(G)\leq\omega(G)-1$ or $G$ is isomorphic to $K_2\square K_{\omega(G)}$.
\end{theorem}

Let $F$ be a graph, and let ${\cal F}$ be the class of $F$-free graphs. By taking ${\cal C}$ to be the (diamond, bull)-free graphs, and applying Theorem~\ref{diamond-1} to ${\cal C}$ and ${\cal F}$, we can immediately obtain the following theorem.

\begin{theorem}\label{diamond-3}
	Let $F$ be a graph and $k\geq2$ be an integer. If $(F,C_3)$-free graph is $k$-colorable, then $\chi(G)\leq \max\{k,\omega(G)\}$ whenever $G$ is ($F$, diamond, bull)-free.
\end{theorem}

In \cite{GHM03}, Gravier {\em et al} proved that every $(P_t,C_3)$-free graph is $(t-2)$-colorable, where $t\geq6$. By Theorem~\ref{diamond-3}, we have that every $(P_t$, diamond, bull)-free graph $G$ satisfies that $\chi(G)\leq\max\{t-2,\omega(G)\}$, where $t\geq6$. This answers a special case of the question of Schiermeyer \cite{S} which stated that whether there exists a constant $c\geq1$ such that every ($P_7$, diamond)-free graph $G$ satisfies that $\chi(G)\leq \omega(G)+c$ or not. This answers a special case of the question of Ju and Huang \cite{JH} stated that for any forest $T$, whether $(T,K_t-e)$-free graph is nearly optimal colorable (this concept has been defined in \cite{JH}, and notice that diamond is $K_4-e$) or not. 

\medskip

\noindent\textbf{{\em Proof of Theorem~$\ref{diamond-2}$}: } Let $G$ be a connected (bull, diamond)-free graph with $\omega(G)\geq3$. Let $\omega(G)=k$, and suppose $\delta(G)\geq k$. Let $A=\{v_1,v_2,\cdots,v_{k}\}$ be a maximum clique of $G$, and let $N_i=\{v\in N(A)~|~N(v)\cap A=\{v_i\}\}$, where $1\leq i\leq k$. We first show that
\begin{equation}\label{e-1}
	N(A)=\cup_{i=1}^{k}N_i.
\end{equation}

It is certain that $\cup_{i=1}^{k}N_i\subseteq N(A)$. Let $x\in N(A)$, and suppose by symmetry that $x\sim v_1$. Since $A$ is a maximum clique of $G$, we may suppose that $x\not\sim v_2$ by symmetry. If $x\sim y$ for some $y\in A\setminus\{v_1\}$, then $\{x,v_1,v_2,y\}$ induces a diamond. So, $x\in N_1$. This proves (\ref{e-1}).

Let $i$ and $j$ be two distinct integers in $\{1, \ldots, k\}$, let $x\in N_i$ and $y\in N_j$. If $x\not\sim y$, then  for some $h\in\{1,\cdots, k\}\setminus\{i,j\}$, $\{x, y, v_i, v_j, v_h\}$ induces a bull. Therefore,
\begin{equation}\label{e-2}
	\mbox{for $1\leq i<j\leq k$, $N_i$ is complete to $N_j$.}
\end{equation}

Since $\delta(G)\geq k$, we have that $N_i\ne\emptyset$ for each $i$ in $\{1, \ldots, k\}$. If $|N_1|\geq2$, let $x, x'\in N_1$, then for $y\in N_2$ and $z\in N_3$, $\{x, x'\}$ is complete to $\{z, y\}$ and $z\sim y$ by (\ref{e-2}).  Now, we have a diamond on $\{x,x',z,y\}$ if $x\not\sim x'$, and have a diamond on $\{x,x',v_1,y\}$ if $x\sim x'$. Therefore, $|N_1|=1$, and similarly, $|N_i|=1$ for $2\leq i\leq k$.

It is certain that $N$ is also a maximum clique of $G$. Therefore, by symmetry and (\ref{e-1}), we have that $N(N)=A$. And so $V(G)=A\cup N$. This completes the proof of Theorem~\ref{diamond-2}. \qed

\medskip

\noindent\textbf{{\em Proof of Theorem~$\ref{diamond-1}$}:} Let $k\geq2$ be an integer. Let ${\cal C}$ be the class of (bull, diamond)-free graphs, and let ${\cal F}$ be a hereditary class of graphs satisfying $\chi(F)\leq k$ if  $F\in {\cal F}$ and $\omega(F)\leq 2$.  Let $G\in {\cal C}\cap {\cal F}$ be a graph. We will prove the theorem by induction on $|V(G)|$. The conclusion is certainly true when $G$ has a small number of vertices.

If $G$ is disconnected, then there exists a component $G'$ of $G$ such that $\chi(G)=\chi(G')$, and by the inductive hypothesis, the theorem holds. So, we may assume that $G$ is connected. By Theorem~\ref{diamond-2}, either $\omega(G)\leq2$, or $\delta(G)\leq \omega(G)-1$, or $G$ is isomorphic to $K_2\square K_{\omega(G)}$.

If $\omega(G)\leq2$ then $\chi(G)\leq k$. If $G$ is isomorphic to $K_2\square K_{\omega(G)}$, then $\chi(G)=\omega(G)$.

So, we suppose that $\delta(G)\leq\omega(G)-1$. Let $v\in V(G)$ with $d(v)\leq \omega(G)-1$. By induction, $\chi(G-v)\leq\max\{k,\omega(G-v)\}$. If $\omega(G-v)\geq k$, then $\chi(G-v)=\omega(G-v)\leq \omega(G)$, and thus we have that $\chi(G)=\omega(G)$ since $d(v)\leq\omega(G)-1$. If $\omega(G-v)<k$, then $\chi(G-v)\leq k$, and thus $\chi(G)\leq k$ since $d(v)\leq\omega(G)-1\leq \omega(G-v)<k$. Therefore, $\chi(G)\leq\max\{k,\omega(G)\}$. This completes the proof of Theorem~\ref{diamond-1}. \qed

\section{Conclusions}

In  studing the relationship between the $\chi$-bounded class of graphs  and the polynomially $\chi$-bounded class of graphs, Chudnovsky {\em et al} \cite{CCDO2023} introduced the concept of Pollyanna, and proved that some special classes of graphs  are Pollyanna (see Theorem~\ref{Chud-pollyanna}). Motivated by these concepts and results, we show that the classes of ($C_4$, bull)-free graphs, ($C_4$, hammer)-free graphs, and (bull, diamond)-free graphs are linear-Pollyanna. Among the problems proposed by Chudnovsky {\em et al} \cite{CCDO2023},  we are particularly interested to the Pollyanna-ness of ($C_4$, diamond)-free graphs and of odd antihole-free graphs.

Notice that the class of ($C_4,C_5$)-free graphs is  a subclass of the odd antihole-free graphs. It may be true that the class of ($C_4$, diamond)-free graphs and the class of ($C_4, C_5$)-free graphs are both linear-Pollyanna. Since the class of bowtie-free graphs is Pollyanna, it would be nice if someone can determine whether the class of ($C_4$, bowtie)-free graphs is linear-Pollyanna or not.

\bigskip

\noindent{\bf Acknowledgement}: We thank the referees for their valuable and helpful comments.

\renewcommand{\baselinestretch}{1.1}

\end{document}